\documentclass[12pt]{amsart}

\usepackage{environ, xcolor}
\NewEnviron{commentA}{}

\newcommand\Aoff{\RenewEnviron{commentA}{}}
\Aoff{} %Change this to Aoff for professional version, Aon to see extra comments

\usepackage{amsmath, amssymb}
\usepackage{array}
\usepackage[frame,cmtip,arrow,matrix,line,graph,curve]{xy}
\usepackage{graphpap, color, paralist, pstricks}
\usepackage[mathscr]{eucal}
\usepackage[pdftex]{graphicx}
\usepackage[pdftex,colorlinks,backref=page,citecolor=blue]{hyperref}
\usepackage{cleveref}
\usepackage{tikz, verbatim}

\newtheorem{theorem}{Theorem}[section]

\newtheorem{lemma}[theorem]{Lemma}

\theoremstyle{definition}
\newtheorem{definition}[theorem]{Definition}
\newtheorem{example}[theorem]{Example}

\newtheorem{question}[theorem]{Question}
\newtheorem{remark}[theorem]{Remark}

\usepackage[margin=1in]{geometry} 
\usepackage{amsmath,amssymb,mathtools, mathrsfs,esint,tikz-cd,verbatim}

\newcommand{\Z}{{\mathbb Z}}

\newcommand{\C}{{\mathbb C}}

\newcommand{\on}[1]{\operatorname{#1}}

\newcommand{\codim}{{\on{codim }}}

\DeclareFontFamily{U}{mathb}{\hyphenchar\font45}
\DeclareFontShape{U}{mathb}{m}{n}{
      <5> <6> <7> <8> <9> <10> gen * mathb
      <10.95> mathb10 <12> <14.4> <17.28> <20.74> <24.88> mathb12
      }{}
\DeclareSymbolFont{mathb}{U}{mathb}{m}{n}
\DeclareMathSymbol{\righttoleftarrow}{3}{mathb}{"FD}
\newcommand{\actsfromleft}{\mathrel{\reflectbox{$\righttoleftarrow$}}}



\title{The dual complex of a $G$-variety}
\author{Louis Esser}

\address{Department of Mathematics, Princeton University, Fine Hall, Washington Road, Princeton, NJ 08544-1000, USA}
\email{esserl@math.princeton.edu}

\begin{document}
\begin{abstract}
We introduce a new invariant of $G$-varieties, the dual
complex, which roughly measures how divisors in the 
complement of the free locus intersect.  We show that the top
homology group of this complex is 
an equivariant birational invariant
of $G$-varieties.  As an application, we demonstrate the 
non-linearizability of certain large abelian group actions
on smooth hypersurfaces in projective space of any dimension and 
degree at least $3$.
\end{abstract}

	\maketitle

\section{Introduction}

Let $G$ be a finite group which acts on a smooth projective variety $X$.
A major problem in algebraic geometry is the classification of such 
$G$-varieties up to $G$-birational equivalence, where
$X$ and $Y$ are \textit{$G$-birationally equivalent} if there are
invariant open sets $U \subset X$, $V \subset Y$ and a
$G$-equivariant isomorphism $\varphi: U \xrightarrow{\cong} V$.

The geometry of $G$-varieties
may be packaged into invariants that
distinguish $G$-birational equivalence classes.
This circle of ideas
begins with the observation that the presence of
an $H$-fixed point is a $G$-birational invariant of $X$
for $H \subset G$ abelian, over an algebraically
closed field of characteristic zero \cite[App. A]{RY1}.
The recent introduction of \textit{equivariant 
Burnside groups} \cite{KT1} has led to 
many new advances in this area.  The $G$-birational invariants
which take values in these groups incorporate 
the birational types of fixed loci of all abelian subgroups
$H \subset G$, actions of $H$ on
the normal spaces to these fixed loci, and 
residual actions on their orbits.  They
have found numerous applications in classifying 
$G$-varieties of low dimension,
such as projective representations of finite groups
\cite{KT2, TYZ2}, algebraic tori \cite{KT3}, 
and actions on other
rationally connected threefolds \cite{CYZ}.

Another tool with significant applications in 
birational geometry is the \textit{dual
complex} of a simple normal crossings divisor. 
This is a simplicial (quasi-)complex which describes 
how the irreducible components of the divisor
intersect.  Such objects may be associated to (log 
resolutions of) singularities \cite{dFKX} and Calabi-Yau pairs
\cite{KX} in such a way that the homotopy type of
the complex is
independent of the resolution chosen. 

Existing $G$-birational invariants
do not directly measure how components
of the non-free locus of a $G$-variety
intersect.  Therefore, in this paper, we introduce
a notion of dual complexes for $G$-varieties,
which roughly records intersection data for the part of the non-free locus
that is an snc divisor.
We show that the top homology group of this complex 
is a $G$-birational invariant (\Cref{biratinvariance})
and use this result to find new examples of
non-birational actions.
In particular, we prove
that actions of rank $n$ abelian subgroups 
$G \subset \mathrm{Aut}(\mathbb{P}^{n+1})$
on invariant smooth hypersurfaces 
$X \subset \mathbb{P}^{n+1}$
of degree $d \geq 3$
are \textit{never} $G$-birational
to a toric action on a smooth toric variety,
apart from some exceptional
actions on cubic surfaces (\Cref{hypersurfaces}).

\subsection{Notation}
$G$ will always denote a finite group. We will only 
consider smooth and projective $G$-varieties
$X$ with faithful action.
All varieties are over an algebraically closed field $k$
of characteristic zero.  For a subgroup $H \subset G$, $X^H$ is
the locus in $X$ on which $H$ acts trivially.
We say that a $G$-variety has \textit{abelian stabilizers}
if for every closed point $x \in X$, 
$\mathrm{stab}_G(x) \coloneqq \{g \in G: g \cdot x = x \}$
is abelian. The \textit{rank} of a finite abelian group
$H$ is the minimal number of generators of $H$.

\noindent \textit{Acknowledgements.} Thank you to Nathan Chen,
J\'{a}nos Koll\'{a}r, Burt Totaro, Ravi Vakil,
and Ziquan Zhuang for useful conversations and comments.
Thank you to Yuri Tschinkel for suggesting an alternative
approach to \Cref{hypersurfaces}.

\section{Dual Complexes}
\label{sec:main}
Suppose $X$ is a smooth projective
$G$-variety with abelian stabilizers.  By Luna's
slice theorem \cite[Lemme III.1]{Luna}, for any closed point
$p \in X$, there is an invariant affine open neighborhood $U$ of
$p$ and a $\mathrm{stab}_G(p)$-equivariant 
morphism $\varphi:U \rightarrow T_p X$
such that $\varphi$ is \'{e}tale at $p$
and $\varphi(p) = 0$.  In
other words, the action is ``\'{e}tale-locally linear".
Therefore, for any subgroup $H \subset \mathrm{stab}_G(p)$,
$T_p X^H \cong (T_p X)^H$ and $X^H$ is
smooth at $p$ (though globally
$X^H$ may be a union of connected components
of different dimensions).  It follows that
the representation of $\mathrm{stab}_G(p)$ on $T_p X$
is faithful.  Thus, for any abelian
group $H \subset G$, any connected component of $X^H$ has
codimension at least $\mathrm{rank}(H)$ in $X$.

The connected components $Z$ of $X^H$ for
which $\codim(Z) = \mathrm{rank}(H)$ are 
particularly well-behaved, by the following lemma.

\begin{lemma}
\label{snc}
Let $X$ be a $G$-variety with abelian stabilizers,
$H \subset G$ an abelian subgroup, and $Z$
a connected component of $X^H$ such that
$\codim(Z) = \mathrm{rank}(H) = k$.
Then $Z$ is a component of an intersection 
$D_1 \cap \cdots \cap D_k$ of divisors, where the general
point $x \in D_i$ has nontrivial
cyclic stabilizer $H_i$, and
$H = \bigoplus_{i = 1}^k H_i$.  
Moreover, these divisors have simple normal
crossings in a neighborhood of $Z$.
\end{lemma}

\begin{proof}
For any point $p \in Z$, we may diagonalize the
action of $H$ on $T_p X$ to 
$T_p Z \oplus \chi_1 \oplus \cdots \oplus \chi_k$, where
the $\chi_i$ are characters of $H$ that do not depend on the 
choice of point $p \in Z$.  The representation
$\chi_1 \oplus \cdots \oplus \chi_k$ of $H$ must
be faithful.  Since $H$ has rank $k$, it follows that
$H = \bigoplus_{i=1}^k H_i$, where
$$H_i \coloneqq \ker(\chi_1) \cap \cdots \cap \widehat{\ker(\chi_i)}
 \cap \cdots \cap \ker(\chi_k)$$
is a nontrivial cyclic group for each $i$.
The subgroup $H_i$ is the stabilizer of a hyperplane
in $T_p X$, so this hyperplane
is the tangent space to a divisor $D_i$ in $X^{H_i}$
passing through $p$.  The divisors $D_1,\ldots,D_k$
have simple normal crossings at $p$ since all of 
their intersections are smooth of the expected 
dimension near $p$ \cite[Lemma 0BIA]{stacks}.
\end{proof}

We refer to the components $Z$ of $X^H$ with 
$\codim(Z) = \mathrm{rank}(H)$
as the \textit{maximal rank strata} of the $G$-variety $X$.
Next, we construct a complex measuring the
intersections of these special strata.

\begin{definition}
The \textit{dual complex} $\mathcal{D}_G(X)$ of a smooth
projective $G$-variety $X$ of dimension $n$ is a CW complex
constructed as follows.  The $0$-skeleton of 
$\mathcal{D}_G(X)$ is a union of points indexed by
divisorial components of $X^H$, as $H$ ranges over all nontrivial cyclic
abelian subgroups of $G$. At the $k$th step, attach a $(k-1)$-simplex for each 
connected component $Z$ of $X^H$ with $\codim(Z) = \mathrm{rank}(H)$,
for $H \subset G$ abelian of rank $k$.
By \Cref{snc}, $Z$ is a connected
component of the intersection $D_1 \cap \cdots \cap D_k$,
where the divisors $D_1,\ldots,D_k$ each have nontrivial
stabilizer and simple normal crossings along $Z$.  The new
simplex is attached along the simplices corresponding to
$D_1,\ldots,D_k$ and the components of their 
intersections containing $Z$.
\end{definition}

The result of this construction is a \textit{simplicial quasicomplex}
$\mathcal{D}_G(X)$ of dimension at most $n-1$.  That is,
it is a CW complex whose cells are simplices such that the
collection of simplices is closed under the face relation and
the intersection of any two is a union
of some of their faces. We refer to 
\cite[Section 7]{ABW} for more about quasicomplexes
and related terminology.
\begin{commentA}
We next review a some notions for quasicomplexes which
naturally arise when considering the dual complexes of 
blowups.

\begin{definition}
Let $\Sigma$ be a simplicial quasicomplex. We use the following 
terminology:
\begin{itemize}
    \item The \textit{star} $\mathrm{Star}(\Sigma,\Delta)$ for
    $\Delta \in \Sigma$ is the set of all faces of $\Sigma$
    containing $\Delta$
    \item The \textit{closure} $\bar{S} \subset \Sigma$
    of a subcollection $S \subset \Sigma$
    consists of every simplex in $S$ and its faces
    \item The \textit{link} of $\Delta$ is 
    $L(\Sigma,\Delta) \coloneqq \overline{\mathrm{Star}(\Sigma,\Delta)} \setminus 
    \mathrm{Star}(\Sigma,\Delta)$
    \item The \textit{stellar subdivision} of a complex $\Sigma$
    in a simplex $\Delta$ with barycenter $p$ is the simplicial subcomplex
    $$\Sigma \setminus \mathrm{Star}(\Sigma,\Delta) \cup
    \{\mathrm{conv}(p,\Delta) | \Delta \in L(\Sigma,\Delta) \}$$
\end{itemize}
\end{definition}
\end{commentA}
\begin{theorem}
\label{biratinvariance}
The homology group $H_{n-1}(\mathcal{D}_G(X))$ is a
$G$-birational invariant of smooth projective
varieties of dimension $n$ with abelian stabilizers.
\end{theorem}

\begin{proof}
By the $G$-equivariant weak factorization theorem \cite{AKMW},
any birational map of smooth projective $G$-varieties
may be decomposed into a sequence of blowups and blowdowns
in smooth centers.  If the varieties both
have abelian stabilizers, we may guarantee
that the intermediate varieties have 
abelian stabilizers as well (cf. \cite[Proposition 3.6]{KT1}).
\begin{commentA}
Indeed, that proposition shows that any two varieties $X$
and $X'$
satisfying Assumption 2 of \cite{KT1} are connected
by a series of blowups and blowdowns in smooth centers,
where each intermediate variety satisfies Assumption 2 
as well.  One of the conditions of Assumption 2 is having
abelian stabilizers.  Though this isn't sufficient to 
satisfy Assumption 2, there are equivariant blowups 
of $X$ and $X'$ satisfying Assumption 2 \cite[Proposition 3.3]{KT1},
and \textit{these} can be connected by intermediate varieties
satisfying the assumption.  Since a blowup of a 
variety with abelian stabilizers has abelian stabilizers, 
this completes the argument.
\end{commentA}
Therefore, it suffices to show that the
homology group $H_{n-1}$ of the dual complex is unchanged 
by a single blowup $X' \coloneqq \mathrm{Bl}_W(X) \rightarrow X$
of a variety $X$ of dimension $n$ with abelian stabilizers
in a $G$-stable center $W$, which is a disjoint union
of smooth irreducible varieties. We may assume that $G$
acts transitively on the connected components of $W$.

There are two cases to consider.  First, suppose
that the rank of the stabilizer of a general point in $W$
equals $\codim(W)$.  Then each component of $W$
is a maximal rank stratum.  We claim that 
in this case $\mathcal{D}_G(X')$ is
the stellar subdivision of the complex
$\mathcal{D}_G(X)$ in the simplices corresponding
to the components of $W$.
Note that distinct components of $W$ are 
disjoint, so no simplex in $\mathcal{D}_G(X)$
contains multiple faces corresponding 
to components of $W$.  It follows that we can
perform all these stellar subdivisions 
simultaneously, and it's enough to consider
the neighborhood of a single component.

To see that the effect of the blowup
is a stellar subdivision, let $Z$ be any 
connected component of $W$,
and $D_1, \ldots, D_k$ the 
divisors from \Cref{snc} containing $Z$.  We continue 
to use the notations $H_i$ and $H$ from the lemma for the groups fixing $D_i$ and $Z$ pointwise, respectively. So each
$D_i$ is fixed by a nontrivial cyclic subgroup $H_i \subset G$
and the sum $H$ of these subgroups has rank $k$.
Denote by $E_Z$ the exceptional divisor over $Z$.
Any maximal rank stratum contained in $Z$ is
a component of an intersection of the 
form $D_1 \cap \cdots \cap D_k \cap D_1' 
\cap \cdots \cap D_{\ell}'$ for some additional
divisors $D_1', \ldots, D_{\ell}'$, where the stabilizer
of a general point on each is a nontrivial
cyclic group.  Furthermore, these divisors 
have simple normal crossings 
in a neighborhood of the
intersection.  It's well known that the dual
complex of the blow up of a stratum in an
snc intersection is the stellar subdivision of the original
complex in the simplex
corresponding to that stratum (see, e.g., \cite[Section 7]{ABW}).
It remains to see that
all the strata of the blowup ``count" 
for our notion of the dual complex,
namely that for each
proper subset $I \subset \{1,\ldots,k\}$,
$E_Z \cap D_I \cap D_1' \cap \cdots \cap D_{\ell}'$
is a maximal rank stratum in $X'$, 
where $D_I = \bigcap_{i \in I} D_i$. 

Indeed,
the stabilizer of a general point $p$ in $E_Z$ is the
intersection of $H \subset \mathrm{GL}((N_{Z/X})_p)$
with the scalar transformations
on this normal space, which must be a nontrivial cyclic
group since $\mathrm{rank}(H) = \dim((N_{Z/X})_p)$.
This group and any $|I| < k$ summands $H_i$
of $H$ generate a rank
$|I|+1$ abelian subgroup.  The intersection
$E_Z \cap D_I \cap D_1' \cap \cdots \cap D_{\ell}'$
is therefore fixed by a group of rank
$|I| + 1 + \ell$, as desired.
Since maximal rank strata
are intersections of maximal rank divisors, no
other new maximal rank strata are introduced by
the blowup. Stellar subdivision induces a homeomorphism
of complexes, so the homology of $\mathcal{D}_G(X')$
is the same as that of $\mathcal{D}_G(X)$ in this case.

The second case is that the
stabilizer of a general point in $W$ has rank smaller
than $\codim(W)$.  For this case, we require the
following lemma:

\begin{lemma}
\label{stable}
Let $H = \Z/d_1 \oplus \cdots \oplus \Z/d_n$
be a rank $n$ abelian group acting diagonally on
$\mathbb{A}^n$.  Then any $H$-stable smooth subvariety
$Y \subset \mathbb{A}^n$ through $0 \in \mathbb{A}^n$
is a coordinate subspace.
\end{lemma}

\begin{proof}
Let $Y = V(I)$ for $I \subset k[x_1,\ldots,x_n]$
a prime ideal.  Since $I$ is $H$-invariant,
we may choose a generating set 
$S \coloneqq \{f_1,\ldots,f_t\}$
for $I$ such that $g \cdot f_i = \chi_i(g) f_i$
for some character $\chi_i$ of $H$.

\begin{commentA}
To prove this, let $T \coloneqq \{h_1,\ldots,h_s\}$ be any 
generating set.  We can create a finite $H$-invariant
generating set $S$ of $I$ by taking the $H$-orbits
of all elements of $T$.
Then $\mathrm{span}_k(S)$ is a finite-dimensional
vector space with $H$-action.  We can diagonalize
the action, so that each basis element is preserved
by $H$ up to scaling.  This basis is also a generating
set for $I$ with the required property.
\end{commentA}

Since $Y$ is smooth and passes through $0$, its 
tangent cone is some linear subspace of $\mathbb{A}^n$.  
It is also
$H$-stable, so we may assume
the initial ideal of $I$ (which cuts out the tangent cone)
is $\mathrm{in}(I) = (x_1,\ldots,x_r)$,
where $\dim(Y) = n-r$.
Among the linear forms in $\mathrm{in}(I)$,
only $x_1,\ldots,x_r$ are preserved
up to scaling by $H$, so the generating set contains
an element of the form $x_1 + h_1$, where $h_1$ only has
terms of degree at least $2$.  Since $x_1$
and $h_1$ must have the same $H$-degree, it follows that
$x_1$ divides $h_1$.  Since $I$ is prime,
$x_1 \in I$.  Similarly, $x_2,\ldots,x_r \in I$.
For dimension reasons, $I = (x_1,\ldots,x_r)$,
so $Y$ is a coordinate plane.
\end{proof}

Returning to the proof of \Cref{biratinvariance},
we claim that the assumption on the stabilizers
of points of $W$ means that $W$
does not contain a maximal rank stratum of dimension $0$.
Suppose by way of contradiction that $W$
does contain such a point $p$ with stabilizer $H$,
with $\mathrm{rank}(H) = n$.
There is a neighborhood $U$ of $p$ and an
$H$-equivariant morphism $\varphi: U \rightarrow T_p X$
which is \'{e}tale at $p$ and sends
$p$ to $0$, by Luna's slice theorem.
Restricting to a smaller open $V \subset U$ containing
$p$ where 
$\varphi$ is \'{e}tale, we have that
$\varphi(V \cap W)$ is an $H$-stable closed subvariety
of an open neighborhood of $0$ in $T_p X \cong \mathbb{A}^n$ 
which passes through $0$ and is smooth there.
The group $H$ is rank $n$ and acts linearly on $T_p X$. 
After diagonalizing this action, 
we may apply \Cref{stable}
to conclude that
$\varphi(V \cap W)$ coincides with a coordinate plane
near $0$ in $T_p X$.  
Every coordinate plane is the image of 
a maximal rank stratum, so $W$ and some stratum
$Z$ coincide \'{e}tale-locally near $x$.  Hence
$W$ is a union of maximal rank strata, a
contradiction.

Since $W$ contains no point with stabilizer of rank 
$n$, it also does not contain any stratum containing
such a point.  Therefore, the subcomplex of $\mathcal{D}_G(X)$
consisting of the closures of all
$(n-1)$-simplices remains the same in
$\mathcal{D}_G(X')$, where each stratum is replaced
by its strict transform, and a collection of strata
intersects in $X$ if and only if it does in $X'$.
There could be new cells in $\mathcal{D}_G(X')$
corresponding to exceptional
divisors over $W$ and their intersections
with existing strata.  However, since every point
of $W$ has stabilizer of rank less than $n$, all new
$k$-cells satisfy $k < n-1$ and don't appear in the
boundary of an $(n-1)$-simplex.
Therefore, the 
boundary map on cells of dimension $n-1$
remains unchanged,
and the homology groups $H_{n-1}(\mathcal{D}_G(X))$
and $H_{n-1}(\mathcal{D}_G(X'))$ are isomorphic.
\end{proof}

The other homology groups of the dual complex
are not $G$-birational invariants.  As a simple
example, suppose the action of $G = \Z/2$ on
a smooth variety $X$ of dimension at
least $2$ has an isolated fixed
point $p$.  Since the codimension of $p$
is greater than $1$, the point is not a maximal
rank stratum and is not accounted for in the dual
complex.  However, blowing up $p$ gives a new fixed
divisor, which will increase the number
of connected components of the dual complex by one.

\section{Applications}

In this section, we consider applications of the
invariant $H_{n-1}(\mathcal{D}_G(X))$ defined in \Cref{sec:main}
to the classification of $G$-actions up
to birational equivalence.
First, note that this group always
vanishes unless the dual complex contains
$(n-1)$-cells, which can only occur if $G$
contains an abelian group of rank at least
$n = \dim(X)$.  Toric varieties form
one natural class with actions of this
kind.

\begin{example}[Smooth Toric Varieties]
\label{toric_ex}
Let $X$ be a smooth projective toric
variety of dimension $n$ and
$G \subset \mathbb{G}_{\on{m}}^n$ a full rank finite
subgroup of the torus.  The non-free locus for the 
$G$-action is precisely the toric boundary, and each
stratum of the boundary has maximal rank by the assumption that
$G$ has rank $n$. 
\begin{commentA}
For instance, consider any toric divisor of $X$, corresponding
to a ray with primitive point $r = (c_1,\ldots,c_n)$ in the lattice
$N$ of one-parameter subgroups of $\mathbb{G}_{\on{m}}^n$.  Our
claim amounts to the fact that any such one-parameter subgroup
intersects $G$ nontrivially.  The group fixing
$D_r$ for $r$ as above is 
$\{(t^{c_1},\ldots,t^{c_n}): t \in \C^*\}$.  Suppose $G$
has elementary abelian subgroup $(\Z/p)^{\oplus n}$.  Then
plugging in $t = \zeta_p$, we find a nontrivial
element $(\zeta_p^{c_1},\ldots,\zeta_p^{c_n})$ in the
intersection of $G$ with the one-parameter subgroup
(it is nontrivial because not all $c_i$ can be divisible by $p$).
The rays spanning any cone in the fan of $X$ form part
of a $\Z$-basis for $N$, so we can use a similar argument
to find a subgroup of the appropriate rank fixing the corresponding
troric stratum. 
\end{commentA}
Since $X$ is smooth and projective, 
its fan is the span of the faces 
of some simplicial polytope $P$ in the lattice $N$ of one-parameter
subgroups \cite[Theorem VII.3.11]{Ewald}.
Therefore, the dual complex $\mathcal{D}_G(X)$
is simply $P$.  In particular,
$\mathcal{D}_G(X)$ is homeomorphic to the sphere $S^{n-1}$,
so the top homology of $\mathcal{D}_G(X)$ 
is $H_{n-1}(S^{n-1}) \cong \Z$.
\end{example}

A faithful $G$-action on a variety $X$ is \textit{linearizable} 
if there is a $G$-birational map $X \dashrightarrow
\mathbb{P}(V)$, where the action on $\mathbb{P}(V)$
is via a faithful linear representation 
$G \rightarrow \mathrm{GL}(V)$.  When $G$ is abelian,
we may diagonalize the $G$-action on $V$ to show that
$G$ is a subgroup of the torus in $\mathbb{P}(V)$.
If $G$ also has rank $n = \dim(X)$, then
\Cref{biratinvariance} and \Cref{toric_ex} show
that $H_{n-1}(\mathcal{D}_G(X)) \cong \Z$ whenever 
the action on $X$ is linearizable.

Rank $n$ abelian subgroups also act on many
hypersurfaces in $\mathbb{P}^{n+1}$.  We
demonstrate that these actions are
never linearizable when the degree
of the hypersurface is at least $3$,
except for the case of cubic surfaces.

\begin{theorem}
\label{hypersurfaces}
Let $n$ and $d$ be positive integers such that
$n \geq 1$, $d \geq 3$, and $(n,d) \neq (2,3)$. Then the action
of any rank $n$ abelian group $G$ 
on a smooth hypersurface $X \subset \mathbb{P}^{n+1}$
of degree $d$ is not
$G$-birational to an action by a subgroup of the torus
on a smooth projective toric variety.  In particular, the action
is not linearizable.
\end{theorem}

\begin{proof}
Smooth plane curves of degree at least $3$
and smooth surfaces of degree at least $4$
in $\mathbb{P}^3$ are not rational, so
we may immediately assume that $n \geq 3$.
In that setting, the Grothendieck-Lefschetz theorem
\cite[Expos\'{e} XII, Corollaire 3.6]{SGA2}
implies that the group $G \subset \mathrm{Aut}(X)$
lifts to a subgroup 
$G \subset \mathrm{Aut}(\mathbb{P}^{n+1}) 
\cong \mathrm{PGL}_{n+2}(k)$.

Every rank $n$ abelian group $G$ contains an elementary
abelian subgroup $(\Z/p)^{\oplus n}$, for some prime $p$.  If $X$
is not $H$-birational to a toric action for $H \subset G$,
it certainly will not be $G$-birational to such an action, 
so we can and do
assume $G \cong (\Z/p)^{\oplus n}$ from now on.

The group $G \subset \mathrm{PGL}(V)$
has some central extension $\tilde{G} \subset \mathrm{GL}(V)$,
where $V \cong k^{n+2}$.
We may assume that there is a $G$-fixed point $p$ on $X$,
or else $H_{n-1}(\mathcal{D}_G(X)) = 0$,
while toric actions
satisfy $H_{n-1} \cong \Z$.  The line $\ell_p \subset k^{n+2}$
is a one-dimensional sub-$\tilde{G}$-representation of $V$. 
After twisting $V$ by this representation 
(which doesn't change the image $G$), 
we can assume that $\tilde{G}$ acts trivially on $\ell_p$.
It follows that $G$ lifts to $\mathrm{GL}_{n+2}(k)$,
and we can diagonalize the action of $G$ on $V$ so that
$V \cong \chi_{\mathrm{triv}} \oplus \chi_1 \oplus 
\cdots \oplus \chi_{n+1}$.
In other words, $\mathbb{P}(V)$ is a compactification
of the linear representation $\chi_1 \oplus \cdots \oplus \chi_{n+1}$
of $G$.

We claim that two of $\chi_{\mathrm{triv}}, \chi_1,\ldots,\chi_{n+1}$
must be the same.
Indeed, if the characters are all distinct, the only $G$-fixed 
points on $\mathbb{P}^{n+1}$ are the coordinate points.
We may assume $X$ contains the coordinate point $p = [1:0:\cdots:0]$
and that $G$ acts on the corresponding line by the trivial character.
Then $T_p X \subset T_p \mathbb{P}^{n+1}$ is a 
sub-$G$-representation.  Since $T_p \mathbb{P}^{n+1}$
is a sum of distinct characters, $T_p X$ coincides with
the tangent space to a coordinate hyperplane,
which we can assume is $\{x_1 = 0\}$ after relabeling.
It also follows that
$T_p X \cong \chi_2 \oplus \cdots \oplus \chi_{n+1}$ is
faithful, so $\chi_2, \ldots, \chi_{n+1}$ form a 
basis for $G^* \cong \mathbb{F}_p^n$.  By \Cref{snc},
$p$ must be the intersection of $n$ divisors on $X$ which
are each fixed by a nontrivial subgroup of $G$. 
These divisors are precisely 
$\{x_2 = 0\} \cap X, \ldots, \{x_{n+1} = 0\} \cap X$.
Indeed, only coordinate strata in $\mathbb{P}^{n+1}$
are fixed by a nontrivial subgroup,
and $X$ cannot contain any codimension $2$
coordinate plane of $\mathbb{P}^{n+1}$ when it is smooth of 
degree $d > 1$ and $n \geq 3$ \cite[Corollary 6.26]{3264}.

Since the hyperplanes $\{x_2 = 0\}, \ldots, \{x_{n+1} = 0\}$ all are
fixed by a nontrivial subgroup, the collection
$\chi_1, \ldots, \widehat{\chi_i},\ldots,\chi_{n+1}$
does \textit{not} form a basis of $G^*$ for $i = 2, \ldots, n+1$.
This implies $\chi_1$ is trivial, contradicting the
assumption that $\chi_{\mathrm{triv}}, 
\chi_1,\ldots,\chi_{n+1}$ are distinct.

We may therefore assume that the $G$-action on $V \cong k^{n+2}$
has the following form: each summand 
of $(\Z/p)^{\oplus n}$ acts by multiplication 
by $p$th roots of unity on a distinct coordinate $x_2, \ldots, x_{n+1}$,
and the action on the first two coordinates is trivial.  Since $X$
cannot contain a codimension $2$ coordinate plane of $\mathbb{P}^{n+1}$,
the only divisors fixed by nontrivial subgroups in $X$
are the intersections
$D_i \coloneqq X \cap \{x_i = 0\}$, $i = 2,\ldots,{n+1}$,
where $D_i$ is fixed by the $i$th summand $\Z/p$ of $G$.
Every intersection of $k < n$ of these divisors
has a rank $k$ stabilizer, so it is smooth and
irreducible of dimension $n - k \geq 1$.  Thus,
the $(n-2)$-skeleton of $\mathcal{D}_G(X)$ is
the same as the $(n-2)$-skeleton of the standard
$(n-1)$-simplex.  However, $D_2 \cap \cdots \cap D_{n+1}$
is a union of $d$ distinct points (or else the divisors do
not intersect transversely), so $\mathcal{D}_G(X)$ has
$d$ top-dimensional simplices with identical boundary.
Therefore, $H_{n-1}(\mathcal{D}_G(X))
\cong \Z^{d-1}$.  By \Cref{biratinvariance} and
\Cref{toric_ex}, $X$ is not $G$-birational to
a smooth toric variety with action by a full rank finite subgroup
of the torus whenever
$d \geq 3$. This
implies in particular that
the $G$-action on $X$ is not linearizable. \end{proof}

\begin{commentA}
\begin{remark}
A related question asks whether any action of $G$ on a hypersurface
$X$ can be \textit{projectively} linearizable, i.e., $G$-equivariantly
isomorphic to any $G$ action on $\mathbb{P}^n$, not necessarily linear.
The techniques of \Cref{hypersurfaces} do not directly answer
this question, but we note that projectively linearizable and 
linearizable are frequently equivalent for a fixed pair $(n,d)$.
In particular, we claim that they are equivalent when $n \geq 3$,
$d \geq 3$, and $\gcd(n+2,d) = 1$.

Indeed, given the finite group $G \subset \mathrm{Aut}(X)$, we
know $G$ lifts to $\mathrm{Aut}(\mathbb{P}^{n+1})$.  But the
group $G$ thus also acts on $K_{\mathbb{P}^4}$
and $K_X$.  These two line bundles are $\mathcal{O}_X(-n-2)$
and $\mathcal{O}_X(d-n-2)$, respectively.  If $\gcd(d,n+2) = 1$,
then some product of tensor powers of these is $\mathcal{O}_X(1)$, so
we get a $G$-action on $H^0(X,\mathcal{O}_X(1))$ inducing
the action on $\mathbb{P}^{n+1}$, meaning that the $G$-action
is linear.  This in turn means that the \textit{Amitsur group}
$\mathrm{Am}(X,G)$ vanishes \cite[Section 6]{BCDP}.  This group
is also an equivariant birational invariant, so if the action
were projectively linearizable, the corresponding action of $G$
on $\mathbb{P}^n$ would also have vanishing Amitsur group.  Hence
it would be linear.  It follows that projective linearizability
implies linearizability in this case.
\end{remark}
\end{commentA}

\begin{example}[Exceptional Cases]
\label{specialcases}
\Cref{hypersurfaces} does not hold for the case of
quadrics, where $d = 2$, and for cubic surfaces, where 
$(n,d) = (2,3)$.  For quadrics, any $G$-action 
with a $G$-fixed point $p$ is linearizable, because the 
projection from $p$ gives a $G$-equivariant 
birational map to projective space.  Of course,
an action by an abelian group $G$ without a fixed point
is not linearizable, since linearizable actions have fixed points.

There is also a counterexample
for smooth cubic surfaces, which are always rational.
Let $X = \{x_0^3 + x_1^3 + x_2^3 + x_3^3 = 0\} \subset \mathbb{P}^3$
be the Fermat cubic surface and $G \subset \mathrm{Aut}(X)$
the normal Klein four-subgroup of the $S_4$ permuting
the four coordinates.  Each involution in $G$ fixes a
line and three isolated points in $X$ (these involutions
are of type 2B in the classification of Dolgachev-Duncan
\cite{DD}).  No nontrivial element of $G$ fixes a positive
genus curve, so $G$ is conjugate to a subgroup of 
$\mathrm{Aut}(\mathbb{P}^1 \times \mathbb{P}^1)$ or 
$\mathrm{Aut}(\mathbb{P}^2)$ in the plane Cremona group
\cite[Theorem 5]{Blanc}.  Since $G$ also fixes the point
$(1:-1:1:-1) \in X$,
the subgroup must be linearizable.
\end{example}

\begin{example}[Rational Hypersurfaces]
\label{rationalexamples}
For every even $n \geq 2$, the dual complex
can be used to identify nonlinearizable actions
on rational cubic hypersurfaces of dimension $n$.
For instance, let $X$ be the Fermat cubic
hypersurface
$$X = \{x_0^3 + \cdots + x_{n+1}^3 = 0\} 
\subset \mathbb{P}^{n+1}.$$
This hypersurface contains a pair of complementary
linear subspaces of dimension $n/2$, so it is rational
\cite[Corollary 1.5.11]{Huybrechts}.
However, there are many actions of the group
$G \coloneqq (\Z/3)^{\oplus n}$ on this hypersurface, such
as the action by $3$rd roots of unity on the first
$n$ coordinates.  By \Cref{hypersurfaces}, this
action is not linearizable.  It is not even 
linearizable in dimension $2$, since the dual complex
is homeomorphic to $S^1 \vee S^1$. In combination with
\Cref{specialcases}, this gives an example
of a single surface with two actions by
rank $2$ abelian groups such that each has a fixed
point, but one is linearizable and one is not.
\end{example}

\begin{remark}
An alternate approach to part of the proof
of  \Cref{hypersurfaces} was later 
suggested by Yuri Tschinkel,
using the formalism of
\textit{incompressible symbols} in equivariant
Burnside groups.
We sketch the argument below in the case of
$G \coloneqq (\Z/d)^{\oplus n}$ acting
on the first $n$ coordinates of a
hypersurface $X$ of dimension
$n$ and degree $d$.

The invariant in the equivariant Burnside group for this
action contains a divisor symbol of the form
$(H,G/H \actsfromleft k(Y), \beta)$. Here
$Y \coloneqq \{x_0 = 0\} \cap X$ is a smooth hypersurface 
of degree $d$ in $\mathbb{P}^n$ with stabilizer $H \coloneqq \Z/d$,
the action $G/H = (\Z/d)^{\oplus (n-1)} \actsfromleft k(Y)$ 
is the residual action of
$G/H$ on the function field of $Y$, and $\beta$ is a character
of $\Z/d$ in the normal space to $Y$.

Such a divisor symbol is \textit{incompressible} roughly if it 
cannot arise as the symbol for an exceptional divisor of
a blowup of a lower dimensional subvariety (for a precise
definition, see \cite[Definition 3.3]{KT2}). Assuming by
induction the nonlinearizability of the residual 
action on $k(Y)$, one can prove the incompressibility of the
corresponding divisor symbol.  It follows that
such a symbol cannot appear in the invariant associated
to a linearizable action of $G$ on $\mathbb{P}^n$.
\end{remark}

When the degree $d$ is small compared to the 
dimension $n$, it is not known whether smooth hypersurfaces
$X_d \subset \mathbb{P}^{n+1}$ are rational.
Even a very general $X$ is only known to be non-rational
when the degree exceeds roughly $\log_2(n)$ \cite{Schreieder}.
In contrast, \Cref{hypersurfaces} shows that
\textit{no} action by a rank $n$ abelian group on a hypersurface
of dimension $n \geq 3$ and degree $d \geq 3$
is linearizable.
This suggests
the following question:

\begin{question}
\label{hypersurface_question}
For which finite groups $G \subset \mathrm{Aut}(\mathbb{P}^{n+1})$
and positive integers $n, d$ is it true that the $G$-action
on any $G$-invariant smooth hypersurface of degree $d$ 
in $\mathbb{P}^{n+1}$ is not linearizable?
\end{question}
A positive answer for $G$ trivial and some $n, d$ is equivalent
to the statement that no smooth hypersurface of degree $d$
in $\mathbb{P}^{n+1}$ is rational.

Besides the groups in \Cref{hypersurfaces}, some other
interesting cases of \Cref{hypersurface_question} in low
dimensions are known.  
For instance, if $G \cong \Z/3$ is the group acting by
third roots of unity on the first coordinate $x_0$ in
$\mathbb{P}^3$, then the induced action
on any $G$-invariant smooth cubic surface $X$ is
not linearizable, even though $X$ is rational.
Indeed, a cyclic group of prime order in the Cremona
group $\mathrm{Cr}(2)$ fixing a curve of positive genus is
not even \textit{stably} $G$-birational to a linear action
on $\mathbb{P}^2$ \cite[Corollary 1.2]{BP}. 
However, \Cref{hypersurface_question} appears to be
open for low rank abelian groups acting on hypersurfaces
of small degree and large dimension.

\end{document}